\documentclass[12pt]{amsart}

\usepackage{amsmath,amssymb,amsfonts}
\usepackage{graphicx}
\begin{document}
\title{Semiclassical asymptotics and entropy}

\author{Tatyana Barron$^1$, Manimugdha Saikia$^2$}

\address{$^{1,2}$Department of Mathematics, University of Western Ontario, London Ontario N6A 5B7, Canada}

\email{$^1$tatyana.barron@uwo.ca, $^2$msaikia@uwo.ca}

\maketitle

\begin{abstract}
We study the entanglement of quantum states associated with submanifolds of Kaehler manifolds. As a motivating example, we discuss the semiclassical asymptotics of entanglement entropy of pure states on the two dimensional sphere with the standard metric.  
\end{abstract}

\section{Preliminaries} 
In this section, we will review some background and facts, and in the next section we will state and discuss our results. 
  
Let $V$ be a finite dimensional Hilbert space of dimension $d$. A vector $v$ in the Hilbert space $V\otimes V$ is {\it separable} (decomposable), if $v=u\otimes w$ for some $u,w\in V$, and {\it entangled} otherwise. This terminology extends to the (pure) states in ${\mathbb{P}}(V)$. The {\it entanglement entropy} of a unit vector $v\in V\otimes V$ is 
\begin{equation}
\label{eentr}
E(v)=-\sum_{j=1}^d\lambda_j\ln \lambda_j
\end{equation}
where $\lambda_1$,...,$\lambda_d$ are  the eigenvalues of ${\mathrm{Tr}}_2(P_v)$, where $P_v$ is the orthogonal projection onto the one-dimensional linear subspace of $V\otimes V$ spanned by $v$ and the partial trace ${\mathrm{Tr}}_2:{\mathrm{End}}(V\otimes V)\to {\mathrm{End}}(V)$ is defined by 
${\mathrm{Tr}}_2(A\otimes B)=B{\mathrm{Tr}}(A)$ for $A,B\in {\mathrm{End}}(V)$ (and extended by linearity). In (\ref{eentr}) the summation is only over $j$ for which $\lambda_j> 0$ (i.e. the convention is 
$0\ln 0=0$).  

The entanglement entropy (\ref{eentr}) is a nonnegative real number between $0$ and $\ln d$. It is $0$ if and only if $v$ is decomposable. We note that 
\begin{equation}
\label{s1inv}
E(v)=E(\gamma v)
\end{equation}
for any $\gamma\in {\mathbb{C}}$ such that $|\gamma|=1$.  {\it Schmidt decomposition:} given a unit vector $v$, there are orthonormal bases 
$\{ e_j\}$ and $\{ f_j\}$ of $V$, and real numbers $\alpha_1\ge ...\ge\alpha_d\ge 0$ such that $v=\sum _{j=1}^d \alpha_je_j\otimes f_j$. 
We note that $\sum\limits_{j=1}^d\alpha_j^2=1$. 
It follows that 
\begin{equation}
\label{eentrs}
E(v)=-\sum_{j=1}^d\alpha_j^2\ln (\alpha_j^2).  
\end{equation}
The numbers $\alpha_j$ are called the Schmidt coefficients of $v$. The number of nonzero Schmidt coefficients of $v$ is the Schmidt rank of $v$. 

As it is well known in quantum information theory, this discussion can be extended to $V_1\otimes V_2$, where $V_1$ and $V_2$ are two different Hilbert spaces, not necessarily finite-dimensional, to the multipartite case, to mixed states, and to many other ways to quantify entanglement. In this note, we explore the simplest possible setup before pursuing more complicated analysis.

\section{Entropy of quantum states}
Our general philosophy is tied with the 
geometry versus analysis perspective. It  emerges in different ways in mathematics, and those often align with physics-driven ideas. For instance, Lagrangian states or Bohr-Sommerfeld states have been a part of the philosophy of geometric quantization for a long time (see e.g. \cite{bpu}). Gelfand-Naimark theorem as well as other reconstruction theorems had a substantial impact in the development of noncommutative geometry and interesting applications to physics (see e.g. \cite{gmt}). 
Shannon entropy (information-theoretic entropy) has been brought into K\"ahler geometry.  Semiclassical asymptotics, with various geometric aspects, have been addressed in \cite{bp, ce, zf}. In \cite{shk} the emphasis is on automorphisms and orbits.

Informally speaking, while discussing the geometry vs. analysis paradigm, by "geometry" we will mean a smooth manifold $M$, possibly with an additional structure (such as a symplectic form, a complex structure, or a Riemannian metric), and by "analysis" we will broadly understand function spaces on $M$, sections of line bundles, operators, norms, estimates, and so on. We will be looking for invariants, or asking to what extent "analysis" determines "geometry" (or whether one can see how the geometric properties of $M$ manifest in "analysis" on $M$). A completely naive example would be to assign to a subset $S$ of $M$ its characteristic function $\chi_S$ and to observe that if $\chi_{S_1}\ne \chi_{S_2}$, then $S_1\ne S_2$.

Let $L\to M$ be a positive holomorphic hermitian line bundle on a compact complex $m$-dimensional manifold $M$ ($m\in {\mathbb{N}} )$. 
It is ample. Assume $L$ is very ample.  
Denote by $\nabla$ the Chern connection on $L$. The $2$-form $\omega=i {\mathrm{curv}}(\nabla)$ is a K\"ahler form on $M$. (Or, alternatively, we could have stated that for an integral K\"ahler manifold $(M,\omega)$ such line bundle exists.) Let $k$ be a positive integer.  For a unit vector $v$ in the finite-dimensional Hilbert space 
\begin{equation}
\label{hilbs}
H^0(M,L^k)\otimes H^0(M,L^k)
\end{equation}
 we can calculate its entanglement entropy $E_k(v)$. 
 The expression for it is the formula for $E(v)$, either (\ref{eentr}) or (\ref{eentrs}). 
 Here $H^0(M,L^k)$ is the space of holomorphic sections of the $k$-th tensor power of $L$, regarded as a (complex) Hilbert space, with the inner product induced by the pointwise hermitian metric on $L$.   
 Let us also denote by $S_k$ the unit sphere in the Hilbert space (\ref{hilbs}): 
 $$
 S_k=\{ v\in    H^0(M,L^k)\otimes H^0(M,L^k)\ | \ ||v||=1\} . 
 $$
 We also note an isomorphism of Hilbert spaces:
 \begin{equation}
 \label{isohs}
 H^0(M,L^k)\otimes H^0(M,L^k)\simeq H^0(M\times M, L^k\boxtimes L^k).
 \end{equation}
 Here $L^k\boxtimes L^k\to M\times M$ is the holomorphic line bundle 
 $$
 \pi_1^* (L^{\otimes k})\otimes 
 \pi_2^* (L^{\otimes k}),
 $$ 
 where $\pi_m:M\times M\to M$, $m\in \{ 1,2\}$, are the projections onto the first and second factor respectively. 
 The hermitian metric on $L$ induces a hermitian metric on $L^k\boxtimes L^k$. 
 For a subset $\Lambda\subset M$, we will denote by $R_k$ the restriction operator defined by 
$$
s\mapsto s\Bigr | _{\Lambda} 
$$
for $s\in H^0(M\times M, L^k\boxtimes L^k)$.  Due to (\ref{isohs}), we will also write $R_k(v)$ for vectors $v$ in (\ref{hilbs}), not to complicate notation.

 In previous work, we showed, in particular:  

\noindent {\bf Theorem 1.} (\cite{bsai})  Let $L\to M$ be a positive holomorphic hermitian line bundle on a compact complex $m$-dimensional manifold $M$. 
Assume $L$ is very ample. Then, as $k\to\infty$, 
$$
\frac{\int_{S_k}E(v)d\mu_k}{\int_{S_k}d\mu_k}=m\ln k-\frac{1}{2}+\ln\beta+\frac{\gamma}{\beta}\frac{1}{k} + O(1/k^2)
$$
where 
$$
\beta =\int_M\frac{c_1(L)^m}{m!}
$$
$$
\gamma =\int_M\Bigl (c_1(L)+\frac{1}{2}c_1(T^{1,0}M)\Bigr )\frac{c_1(L)^{m-1}}{(m-1)!}
$$
and $d\mu_k$ is the spherical measure on $S_k$. 

Here $c_1(L)$ denotes the first Chern class of $L$ and $c_1(T^{1,0}M)$ denotes the first Chern class of the holomorphic tangent bundle of $M$.  

\

Now, we will aspire to assign to a subset of $M\times M$ (say, $\Lambda$), some analytic construct built from the Hilbert space (\ref{hilbs}) via 
the quantum information concepts, such as entanglement entropy, negativity or entanglement of formation. As a start, in this note, we will consider 
$M={\mathbb{CP}}^1$ with the Fubini-Study metric, $L$ the hyperplane bundle, $\Lambda=S^1$ embedded antidiagonally as specified below, and we will prove the theorem below, with the intent to generalize this statement later to other $M$ and $\Lambda$ and to use Theorem 2 and its proof as a guiding example for the general case.    

\

\noindent {\bf Theorem 2.} Let $M={\mathbb{CP}}^1$ equipped with the Fubini-Study metric and let $L\to M$ be the hyperplane bundle, with the standard hermitian metric. Let $\Lambda\subset  (M-\{ [1:0]\})\times (M-\{ [1:0]\})$ be defined by 
$$
\Lambda=\{ (z,w)\in {\mathbb{C}}\times {\mathbb{C}} \ | \ z=e^{it}; w=e^{-it}; 0\le t\le 2\pi\} ,
$$
where 
$$
M-\{ [1:0]\}\simeq {\mathbb{C}}
$$ 
is the affine chart 
$$
\{ [z_0:z_1]\in {\mathbb{CP}}^1 \ | \ z_1\ne 0 \} 
$$
with the affine coordinate $z=\frac{z_0}{z_1}$. Let $\{ e_j\}$ be the standard orthonormal basis in $H^0(M,L^k)$:
$$
e_j=\sqrt{\frac{(k+1)!}{j!(k-j)!} }z^j
$$
$j\in \{ 0,1,2,...,k\}$. Let for each $k\in{\mathbb{N}}$, 
$$
W_k=\ker (R_k)\cap {\mathrm{span}} \{ e_j\otimes e_j| \ j=0,1,...,k\} .
$$ 
Then

\noindent (a) $\ker (R_1)$ is the span of $\frac{1}{\sqrt{2}}(e_0\otimes e_0-e_1\otimes e_1)$

\noindent (b) The sequence of vectors 
$$
b_k=\frac{1}{\sqrt{1+k^2}}(e_1\otimes e_1-ke_0\otimes e_0)
$$
is in $\ker R_k$ for each $k=1,2,3,...$, and their entanglement entropy is 
\begin{equation}
\label{eentrbk}
E_k(b_k)=-\frac{1}{1+k^2}\ln \frac{1}{1+k^2}-\frac{k^2}{1+k^2}\ln \frac{k^2}{1+k^2}.
\end{equation}
\noindent (c) The sequence of vectors 
$$
c_k=\frac{1}{\sqrt{2}}(e_0\otimes e_0-e_k\otimes e_k)
$$
is in $\ker R_k$ for each $k=1,2,3,...$, and their entanglement entropy is $E_k(c_k)=\ln 2$ for every $k$. 

\noindent (d) The linear subspace $W_k$ is $k$-dimensional. For all odd $k$, 
there is a vector of maximum Schmidt rank in $W_k$ such that its entanglement entropy is $\ln (k+1)$.  For each even $k$, $W_k$ contains a vector whose entanglement entropy is $\ln k$. 

\noindent (e) The endomorphism of the Hilbert space $H^0(M\times M, L\boxtimes L)$ defined by the orthogonal projection onto $\ker (R_1)$ is 
the Berezin-Toeplitz operator $T_f^{(1)}$ with the symbol $f\in C^{\infty}(M\times M)$ given by 
$$
f(z_0, z_1;w_0, w_1)=\frac{9}{2}\frac{(z_0w_0-z_1w_1)(\bar z_0 \bar w_0-\bar z_1\bar w_1)}{(|z_0|^2+|z_1|^2)(|w_0|^2+|w_1|^2)}-2. 
$$ 
{\it Remark.} For $k\in{\mathbb{N}}$ and a smooth function $F: M\times M\to {\mathbb {C}}$,  the Berezin-Toeplitz operator $T_F^{(k)}$ with the symbol $F$, is an endomorphism of the Hilbert space $H^0(M\times M, L^k\boxtimes L^k)$ defined by  
$$
T_F^{(k)}: s\mapsto \Pi_k(Fs)
$$
where $\Pi_k:L^2(M\times M, L^k\boxtimes L^k)\to H^0(M\times M, L^k\boxtimes L^k)$ is the orthogonal projection from the space of $L^2$ sections of 
the line bundle $L^k\boxtimes L^k$ onto the closed subspace of holomorphic sections.  
 \begin{figure}
\begin{center}
\includegraphics[width=18pc]{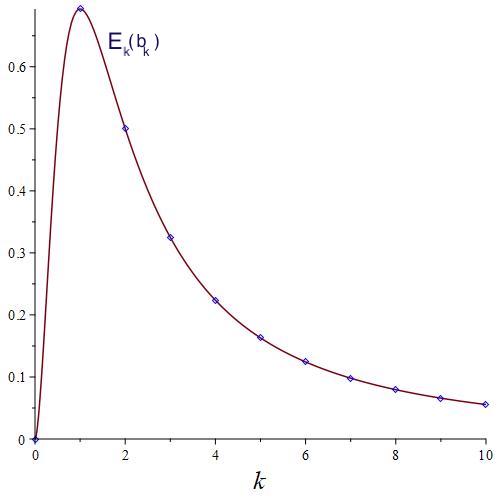}
\end{center}
\caption{The values of (\ref{eentrbk}) for $1\le k\le 10$. }
\end{figure}

\noindent {\it Comments.} 
\begin{description}
\item In Theorem 2(b), the sequence of vectors $b_k$ is asymptotic to $-e_0\otimes e_0$. Accordingly, $E_k(b_k)\to 0$ as $k\to \infty$. 

\item In Theorem 2(c), the vector $c_k$ represents a Bell state for each $k$. 

\item Part (d) shows that the value of $E_k$ on $W_k$ reaches the maximum possible value. Generally, it does not have to be the case (for example,  
the value of entanglement entropy on a $1$-dimensional subspace spanned by a decomposable vector, is zero, or, as another example, 
if we consider the $2$-dimensional subspace spanned by $e_0\otimes e_0$ and $e_0\otimes e_1$, every vector in this subspace is decomposable and its entanglement entropy is zero). 

\item If we realize the linear operators that are relevant to the quantum information theory, as Toeplitz operators, then we would be able to use the asymptotic ($k\to\infty$) results about the spectrum of Toeplitz operators to make conclusions about the semiclassical asymptotics of entanglement.  Part (e) is a demonstration of the first part of this statement in the setting of Theorem 2.  Careful asymptotic analysis in a general situation is upcoming in our  work.  
\end{description}

\section{Proof of Theorem 2.}

\noindent Proof of (a). We observe: 
$$
a \ e_0\otimes e_0+b \ e_0\otimes e_1+c \ e_1\otimes e_0+d \ e_1\otimes e_1
$$
is in $\ker (R_1)$ if and only if $a=d$ and $b=c=0$. The statement follows. 

\noindent Proof of (b) and (c). We observe that, restricted to $\Lambda$, 
 \begin{equation}
 \label{onlambda}
 e_j\otimes e_j={k \choose j}e_0\otimes e_0
 \end{equation}
 for all $0\le j\le k$. Suppose a vector $v$ in (\ref{hilbs}) 
 is of the form 
 $$
 v=\sum_j a_je_j\otimes e_j.
 $$
 Then, the condition for $v$ 
 to be in $\ker (R_k)$ is 
 $$
 \sum_j {k\choose j} a_j=0.
 $$
 We see that the vectors $b_k$ and $c_k$ satisfy this condition and therefore are in $\ker (R_k)$ for all $k$. The values of their entanglement entropy are obtained from their Schmidt coefficients via  (\ref{eentrs}).

 \noindent Proof of (d). The first claim follows from (\ref{onlambda}) and the fact that $e_i\otimes e_j\notin W_k$ for $i\ne j$. 
Now, suppose $v=\sum_{j=0}^k a_je_j\otimes e_j$  is of norm $1$, i.e. 
\begin{equation}
\label{norm1}
\sum_{j=0}^k|a_j|^2=1. 
\end{equation}
Suppose $v$ is in $W_k$ and it is of maximal Schmidt rank. It is equivalent to 
\begin{equation}
\label{inkernel}
\sum_{j=0}^k{k\choose j} a_j=0
\end{equation}
(because it is in $\ker (R_k)$), and $a_j\ne 0$ for all $j$. By a direct calculation, we conclude that the entanglement entropy 
\begin{equation}
\label{entrv}
E_k(v)=-\sum_{j=0}^k|a_j|^2\ln(|a_j|^2).
\end{equation}
Write $a_j=x_j+iy_j$ for $1\le j\le k$. The right hand side of (\ref{entrv}) is a function of  $2k$ real variables (using (\ref{norm1}) ) :
$$
f(x,y)=-(1-\sum_ {j=1} ^k (x_j^2+y_j^2) )\ln (1-\sum_{j=1}^k (x_j^2+y_j^2) )- \sum_{j=1}^k(x_j^2+y_j^2)\ln(x_j^2+y_j^2).
$$
It is $S^1$-invariant with respect to the circle action on (\ref{hilbs}) (see (\ref{s1inv})). To look for the critical points, we consider the equation  
$$
\nabla f=0
$$ 
that leads to 
$$
|a_1|^2=...=|a_k|^2=1-\sum_{j=1}^k|a_j|^2.
$$
When $k$ is odd, there is a solution to these equations that also satisfies (\ref{inkernel}), with $|a_j|^2=\frac{1}{k+1}$ for all $1\le j\le k$ and such that 
$$
a_j=-a_{k-j}
$$
for all $0\le j\le k$. 
When $k$ is even, we can set $a_{k/2}=0$ and choose $a_j$  for all $j\ne \frac{k}{2}$ so that $|a_j|^2=\frac{1}{k}$ and $a_j=-a_{k-j}$. 
The value of $E_k$ is now obtained  from (\ref{entrv}). 

\noindent Proof of (e). The conclusion is obtained by a direct calculation: we apply the matrix of the orthogonal projection onto $\ker(R_1)$ and the linear operator $T_f^{(1)}$ to the four basis vectors. 

\

{\bf Acknowledgements.} Both authors are thankful to the conference organizers for their efforts. The referee's suggestions were helpful and improved the exposition of this paper.

\end{document}